\documentstyle{article}
\begin{document}

\def\be{\begin{equation}}
\def\ee{\end{equation}}
\def\bc{\vskip3mm\begin{center}}
\def\ec{\end{center}\vskip3mm}

\author{Iskander A. Taimanov}
\title{Surfaces of revolution in terms of solitons}

\date{}
\maketitle

\bc
{\bf \S 1. Introduction}
\ec         

In the present article we examine in details global deformations of surfaces 
of revolution via the modified Korteweg--de Vries (mKdV) equations and
the first integrals, of these deformations, regarded as invariants of
surfaces. It is a sequel to
our paper \cite{T} where the general case of modified Novikov--Veselov (mNV)
deformations is considered. Since the main problems are still open, we
show how to solve part of them for surfaces of revolution regarded 
as a toy model for the general case. 

The main tool is the modified Novikov--Veselov equation
\be
U_t = (U_{zzz} + 3U_zV + \frac{3}{2}UV_z) +
(U_{{\bar z}{\bar z}{\bar z}} + 3U_{\bar z}{\bar V} + 
\frac{3}{2}U{\bar V}_{\bar z})
\label{1.1}
\ee
where $V_{\bar z} = (U^2)_z$ and $z = x + \sqrt{-1}y \in {\bf C}$. 
It was introduced by Bogdanov (\cite{Bogd})
and it is a modification of the 
Novikov--Veselov equation (\cite{VN}) in the same
manner as the mKdV equation is a modification of the Korteweg--de Vries
(KdV) equation. It is easy to see that (\ref{1.1}) reduces to the mKdV 
equation for $U(z) = U(x)$. 

In \cite{Kon} Konopelchenko did an important observation that
if a surface is represented locally by the generalised Weierstrass 
representation
\footnote[1]
{This is the representation (\ref{2.13}) for a zero eigenfunction $\psi$ of 
the operator (\ref{1.1}).}
from the operator
\be
L = 
\left(
\begin{array}{cc}
\partial & -U \\
U & {\bar \partial}
\end{array}
\right)
\label{1.2}
\ee
and we take the mNV deformation of the potential 
$U$, then this deformation induces a local deformation of a surface.

In fact, this local representation of surfaces immersed into ${\bf R}^3$
is not new and it seems that it was Eisenhart who first mentioned it
(a short overview is contained in \cite{T}).
Now we know that these formulae give a local representation of arbitrary 
surface.

In \cite{T} we proposed to study the mNV deformations of surfaces globally.
This led to an examination of first integrals of the mNV equations as geometric
quantities. The first of them has the same density as 
the Willmore functional. We correctly defined the mNV flow for
immersed tori and proved that these tori are preserved by this deformation.
This implies that the Willmore functional is the first integral of the 
global soliton deformation of immersed tori.

The interesting feature of the generalised Weierstrass representation that
globalising it we naturally arrive at regarding generic immersed surfaces via
Dirac operators. In fact, global immersions of a surface $\Sigma$ is in 
one-to-one correspondence with zero eigenfunctions of the Dirac 
operator (\ref{1.1}) acting on the spinor bundle to $\Sigma$.
This was shown in \cite{T} (see Propositions 1--4).
\footnote[2]
{Recently this has been rewritten in \cite{KS} in terms of Sullivan's spinor 
representation of minimal surfaces. However for minimal surfaces 
this is the classical Weierstrass representation, taken in the spinor form
it leads to many interesting results (\cite{KS}).}

This observation enables us to introduce a notion of the potential of
an immersed surface. This is the potential $U$ of the Dirac operator.

Recently it was proposed to change a complex function theory
to a quaternionic function theory for regarding generic surfaces
(\cite{KPP}). However this approach is in 
progress now,
it already has led to an interesting explicit relation between Bonnet 
pairs and isothermic surfaces.

In particular, in the present article we prove that all mKdV 
flows preserve tori of 
revolution (Theorem 1), describe its stationary points (Propositions 5 and 6),
describe the deformations of the potentials of surfaces of revolution under
M\"obius transformation (Propositions 7 and 8) and 
a passing to the dual surface (Proposition 9), 
and consider a lot of examples seeking for hints how to generalise 
all reasonings for generic surfaces.   
 
In \S 6 we discuss two interesting 
problems on invariance of 
the first integrals of the mKdV flows taken as 
invariants of surfaces of revolution under the M\"obius transformations 
with centres on the line $Z_2 = 0$ and
a passing to the dual curve in the sense of isothermic surfaces.
We present results of numeric experiments confirming this. 

\bc
{\bf \S 2. Representation of surfaces of revolution via a Dirac operator}
\ec

Let $\Sigma$ be a surface of revolution in the three-dimensional
Euclidean space ${\bf R}^3$ with coordinates $Z_1, Z_2$, and $Z_3$. 

Without loss of generality, assume that $\Sigma$ is obtained by
revolving a curve $\gamma$ lying in the plane $Z_2 = 0$ around the line
$Z_1 = Z_2 =0$. Introduce the conformal parameters $x$ and $y$ on $\Sigma$
as follows 
\be
Z_1 = \theta(x) \cos y, \ \
Z_2 = \theta(x) \sin y, \ \ 
Z_3 = \varphi(x),
\label{2.1}
\ee
where
\be
\theta^2(x) = \left( \frac{d\theta(x)}{dx} \right)^2 +
\left( \frac{d\varphi(x)}{dx} \right)^2.
\label{2.2}
\ee

In these coordinates the first fundamental form is written as
\be
I = \theta^2(x) ( dx^2 + dy^2) = 
\left(\left(\frac{d\theta(x)}{dx}\right)^2 + 
\left(\frac{d\varphi(x)}{dx}\right)^2\right) 
\left( dx^2 + dy^2 \right).
\label{2.3}
\ee

It follows from (\ref{2.2}) that we may correctly define the smooth
functions $\sqrt{\theta-\theta_x}$, $\sqrt{\theta+\theta_x}$, and
$\sqrt{\theta^2-\theta^2_x}$ meaning by them their branches smoothly extended
through zeroes.
By (\ref{2.2}), we mean the function $\varphi_x$ by 
$\sqrt{\theta^2 - \theta^2_x}$ and suppose that
$\sqrt{\theta - \theta_x} \cdot \sqrt{\theta+\theta_x} = \varphi_x$. 

Now, define the following function
\be
U(x) = \frac{\theta - \theta_{xx}}{4\sqrt{\theta^2 - \theta^2_x}}.
\label{2.4}
\ee

We also introduce the vector function 
${\psi}(x,y) = (\psi_1(x,y), \psi_2(x,y))$ as follows
\be
\psi_1(x,y) = r(x) \exp\frac{\sqrt{-1}y}{2}, \ \ 
\psi_2(x,y) = s(x) \exp\frac{\sqrt{-1}y}{2}.
\label{2.5}
\ee
where
\be
r(x) = \sqrt{\frac{\theta - \theta_x}{2}}, \ \ 
s(x) = \sqrt{\frac{\theta + \theta_x}{2}}.
\label{2.6}
\ee

It is checked by straightforward computations that
\be
r_x = -\frac{1}{2}r + 2Us, \ \
s_x = \frac{1}{2}s - 2Ur.
\label{2.7}
\ee
The latter equations are nothing else but the linear problem
\be
{\tilde L}{\tilde \psi} = 0
\label{2.8}
\ee
where
\be
{\tilde L} = \frac{d}{dx} -
\frac{1}{2}
\left(\begin{array}{cc}
-1 & 4U \\
-4U & 1
\end{array}\right)
\label{2.9}
\ee
and
\be
{\tilde \psi} = (r,s).
\label{2.10}
\ee

The linear problem (\ref{2.8}--\ref{2.10}) is rewritten as
\be
L \psi = 0
\label{2.11}
\ee
where
\be
L = 
\left(\begin{array}{cc}
\partial_z & 0 \\
0 & \partial_{\bar z}
\end{array}\right) +
\left(
\begin{array}{cc}
0 & - U \\
U & 0
\end{array}
\right)
\label{2.12}
\ee
and $z = x + \sqrt{-1} y$.

Moreover, the surface $\Sigma$ is represented via the 
eigenfunction :
$$
Z_1 = 
\frac{\sqrt{-1}}{2}
\int \left(({\bar \psi}^2_1 + \psi^2_2) dz  - 
({\bar \psi}^2_2 + \psi^2_1) d{\bar z} \right),
$$
\be
Z_2 = \frac{1}{2}
\int \left(({\bar \psi}^2_1 - \psi^2_2) dz -
({\bar \psi}^2_2 - \psi^2_1) d{\bar z} \right),
\label{2.13}
\ee
$$
Z_3  =  - \int ({\bar\psi}_1\psi_2 dz + \psi_1{\bar\psi}_2 d{\bar z}).
$$

Hence, we arrive at the following conclusion.

{\bf Proposition 1.}
{\sl Every smooth surface of revolution $\Sigma$ is represented via the
zero eigenfunction $\psi$ of the Dirac operator (\ref{2.12}) by
the formulae (\ref{2.13}). Moreover, there exists an explicit 
procedure of constructing this operator by 
(\ref{2.1}--\ref{2.4}).} 
 
Consider an explicit construction of a surface of revolution from
a metric tensor $\theta(x)$ only. We have

{\bf Proposition 2.}
{\sl The following formula
\be
\theta(x) = \exp\int \tau(x') dx'
\label{2.14}
\ee
together with (\ref{2.5}), (\ref{2.6}), and (\ref{2.13}) give 
an explicit procedure of constructing a surface of revolution
from an arbitrary  function 
$\tau(x)$ meeting the condition
\be
\tau^2(x) \leq 1.
\label{2.15}
\ee
If $\theta(x)$ is periodic, $\theta(x) = \theta(x+T)$, then 
the surface of revolution is closed if and only if}
\be
\int^T_0 \sqrt{\theta^2 - \theta^2_x} dx = 0.
\label{2.16}
\ee

For completeness of explanation give explicit formulae for the curvatures
of the surface. 

We have the following formulae for its normal vector 
\be
{\vec N} = \frac{{\vec Z}_x \times {\vec Z}_y}{|{\vec Z}_x \times {\vec Z}_y|} 
= \left(
-\frac{\varphi_x}{\theta}\cos y, - \frac{\varphi_x}{\theta}\sin y,
\frac{\theta_x}{\theta} \right),
\label{2.17}
\ee
its second fundamental form
\be
II = ({\vec N}, d^2{\vec Z}) = 
\left(\frac{\theta\varphi_{xx}}{\theta_x} - 
\varphi_x \right) dx^2 + \varphi_x dy^2,
\label{2.18}
\ee
its Gaussian curvature
\be
K = \frac{\theta^2_x - \theta \theta_{xx}}{\theta^4} = 
- \frac{(\log \theta)_{xx}}{\theta^2},
\label{2.19}
\ee
and its mean curvature
\be
H = \frac{\varphi_{xx}}{2\theta \theta_x} = \frac{2U}{\theta}.
\label{2.20}
\ee

{\bf Examples.}

1) The cylinder of radius $1$ :
$$
\theta(x) = 1, \ \ \varphi(x) = x, \ \ U(x) = \frac{1}{4}.
$$  

2) The unit sphere $Z^2_1 + Z^2_2 + Z^2_3 =1$ :
\be
\theta(x) = \frac{1}{\cosh x}, \ \ 
\varphi(x) = \tanh x, \ \ 
U(x) = \frac{1}{2\cosh x}.
\label{2.21}
\ee

3) The Clifford torus (the revolving curve $\gamma$ is the circle of radius
$1$ with the centre at $(\sqrt{2},0)$):
$$
\theta(x) = \frac{1}{\sqrt{2} - \sin x}, \ \
\varphi(x) = \frac{\cos x}{\sin x - \sqrt{2}}, \ \
U(x) = \frac{\sin x}{2\sqrt{2}(\sqrt{2} - \sin x)}.
$$
 
4) The round torus $T^2_R$
obtained by revolving the circle of radius $1$ with
the centre at $(R,0)$ (this is a generalization of the Clifford torus) :
$$
\theta(x) = R - \sin f(x), \ \
\varphi(x) = \cos f(x), \ \
$$
where 
$$
\frac{df}{R - \sin f} = dx.
$$
In this case the potential
$$
U(x) = \frac{R - 2 \sin f(x)}{4}
$$
satisfies the following equality
\be
U^2_x = \frac{1}{4}\left(2U+\frac{R}{2}\right)^2
\left(1-\left(2U-\frac{R}{2}\right)^2\right).
\label{2.22}
\ee

\bc
{\bf \S 3. The hierarchy of modified Korteweg--de Vries (mKdV) equations}
\ec

We remind some facts about the hierarchy of mKdV equations
(here we follow to formulae from \cite{Schief}).

We consider the periodic problem for these equations and 
denote by $T$ a period.

This hierarchy is defined as
the compatibility conditions for the following linear problems
\be
L(\lambda) \psi = 
\left[
\frac{\partial}{\partial x} - 
\frac{1}{2}
\left(
\begin{array}{cc}
\lambda & q \\
-q & -\lambda
\end{array}
\right)
\right] \psi = 0
\label{3.1}
\ee
and
\be
\left[
\frac{\partial}{\partial t} - 
\frac{1}{2}
K_{n}(\lambda)
\right]\psi = 
\left[
\frac{\partial}{\partial t} -
\frac{1}{2}
\left(
\begin{array}{cc}
A^{(n)} & B^{(n)} \\
C^{(n)} & -A^{(n)}
\end{array}
\right)
\right] \psi = 0.
\label{3.2}
\ee

The entries of $K_n(\lambda)$ are as follows
$$
A^{(n)} = \sum_{k=0}^n A^{(n)}_{2k+1} \lambda^{2k+1},
$$
$$
B^{(n)} + C^{(n)} = 2 \sum_{k=0}^{n-1} S^{(n)}_{2k+1} \lambda^{2k+1}, \ \
B^{(n)} - C^{(n)} = 2 \sum_{k=0}^n T^{(n)}_{2k} \lambda^{2k},
$$
\be
A^{(n)}_{2k+1} = \partial_x^{-1} ( q D^{n-k-1}q_x), \ \ A_{2n+1}^{(n)} = 1,
\label{3.3}
\ee
$$
S^{(n)}_{2k+1} = D^{n-k-1}q_x, \ \ 
T^{(n)}_{2k} = \partial_x^{-1} D^{n-k}q_x ,
$$
where the recursion operator $D$ takes the form
\be 
D = \partial^2_x + q^2 + q_x\partial^{-1}_xq.
\label{3.4}
\ee

The compatibility condition for (\ref{3.1}) and (\ref{3.2}) is
\be
\left[ L(\lambda), \frac{\partial}{\partial t} - 
\frac{1}{2} K_n(\lambda) \right] = 0. 
\label{3.5}
\ee
The $n$-th equation of this hierarchy, the commutation relation (\ref{3.5}),
is written as
\be
\frac{\partial q}{\partial t_n} = D^n q_x.
\label{3.6}
\ee
This representation is not rather correct because it contains the many-valued
operation $\partial^{-1}_x$ but (\ref{3.6}) defines the $n$-th mKdV
equation up to 
a linear combination of the first $n-1$ mKdV equations.

The first equations  are
\be
n = 1\ : \ \ \ q_t = q_{xxx} + \frac{3}{2}q^2q_x,
\label{3.7}
\ee
\be
n = 2\ : \ \ \ q_t = q_{xxxxx}  + \frac{5}{2}q^2q_{xxx} + 10qq_xq_{xx}
+ \frac{5}{2}q^3_x + \frac{15}{8}q^4q_x.
\label{3.8}
\ee

Rewrite the linear problem (\ref{3.1}) for the functions
$$
\chi_1 = \frac{\psi_1 + \sqrt{-1}\psi_2}{\sqrt{2}}, \ \
\chi_2 = \frac{\sqrt{-1}\psi_1+\psi_2}{\sqrt{2}}
$$
and obtain
\be
\left[
\partial_x - \frac{1}{2}\Omega
\left( 
\begin{array}{cc}
\lambda & q \\
-q & -\lambda
\end{array}\right) \Omega^{-1} 
\right] \chi = 
\left[
\partial_x - \frac{\sqrt{-1}}{2}
\left(\begin{array}{cc}
-q & -\lambda \\
\lambda & q 
\end{array}\right)
\right] \chi = 0
\label{3.9}
\ee
where
$$
\Omega = \frac{1}{\sqrt{2}}
\left(\begin{array}{cc}
1 & \sqrt{-1} \\
\sqrt{-1} & 1
\end{array}\right).
$$
Denote $\chi_2$ by $\xi$ and denote $-\lambda^2/4$ by $\mu$.
Now rewrite 
(\ref{3.9}) as
\be
\left( \partial_x + \frac{\sqrt{-1}q}{2} \right)
\left(\partial_x - \frac{\sqrt{-1}q}{2} \right) \xi = -\mu \xi.
\label{3.10}
\ee
Introducing the potential
\be
u = \frac{q^2}{4} - {\sqrt{-1}}\frac{q_x}{2}
\label{3.11}
\ee
we infer that (\ref{3.10}) is the eigenvalue problem for 
the one-dimensional Schr\"odinger operator with this potential
\be
(\partial^2_x + (u+\mu)) \xi = 0.
\label{3.12}
\ee
The transform $q \rightarrow u$ given by (\ref{3.11}) is the Miura transform
assigning solutions to the Korteweg--de Vries (KdV) equation 
\be
u_t = u_{xxx} + \frac{3}{2}uu_x
\label{3.13}
\ee
to solutions to the mKdV equation. 

{\bf Theorem of Miura.}
{\sl Let a function $q(x)$ satisfy (\ref{3.7}).
Then the function $u(x)$ (\ref{3.11}) satisfies (\ref{3.13}).} 

There is a procedure for constructing an infinite family of first integrals of
the KdV equation, 
the Kruskal integrals. It is as follows.
Consider the asymptotic behaviour of $\xi(x,\sqrt{\mu})$ as 
$\mu \rightarrow \infty$:
\be
\xi(x,\sqrt{\mu}) = \exp{(\sqrt{-1}\sqrt{\mu}x + \Phi)}
\label{3.14}
\ee
and rewrite (\ref{3.12}) in terms of (\ref{3.14})
\be 
\Phi_{xx} + \Phi^2_x + u + 2\sqrt{-1}\sqrt{\mu}\Phi_x = 0.
\label{3.15}
\ee
A solution to (\ref{3.15}) is given in the following form
\be
\Phi_x = \sum_{k=1}^{\infty} \frac{R_k}{(2\sqrt{-1}\sqrt{\mu})^k}
\label{3.16}
\ee
where $R_k$ are polynomials in $u$ and its derivatives and are obtained
by the following recursion relations
\be
R_1 = - u, \ \ 
R_{n+1} = -R_{nx} - \sum_{k=1}^{n-1} R_k R_{n-k}.
\label{3.17}
\ee
Since $R_{2n}$ are the derivatives of periodic functions, only the integrals
\be
H_n = -\int_0^T R_{2n+1}
\label{3.18}
\ee 
do not vanish identically.
These are the Kruskal integrals, and the first of them are
$$
H_0 = \int u dx, \ \ H_1 = \int u^2 dx, \ \
H_2  = \int (2u^3 - u^2_x) dx.
$$

Now, we conclude that, by the Miura theorem, these integrals 
generate an infinite family of first integrals of the mKdV equation
\be
{\tilde H}_n (q) = H_n(\frac{q^2}{4}-\frac{\sqrt{-1}q_x}{2}).
\label{3.19}
\ee
which starts with
\be 
{\tilde H}_0 = \frac{1}{4} \int q^2 dx, \ \ 
{\tilde H}_1 = \frac{1}{16} \int (q^4 - 4 q^2_x) dx, \ \ 
\label{3.20}
\ee
$$
{\tilde H}_2 = \frac{1}{32} \int (q^6  - 20 q^2q_x^2 + 8 q_{xx}^2) dx.
$$

Since the Kruskal integrals for the KdV equation are the Hamiltonian functions
for the higher KdV equations and the Miura transform, in fact, assigns to a
solution to the $n$-th mKdV equation a solution to the
$n$-th KdV equation, it is proved that for every $n \geq 1$ 
the integral ${\tilde H}_n$ is the Hamiltonian function 
for the $k$-th mKdV equation. Moreover, since the mKdV-flows commute, 
${\tilde H}_k$ are the first integrals for all equations 
of the mKdV hierarchy. 

We will not dwell on the derivations of these facts because they are standard
and well-known. Probably, this is not explained in details in the main 
textbooks on soliton theory but by using the Miura transform one can easily 
derive this from the analogical results for the KdV hierarchy. For the 
latter we refer, for instance, to \cite{New}.

We summarise the necessary information in the following 
proposition.

{\bf Proposition 3.}
{\sl The functionals $H_k$ given by (\ref{3.15}--\ref{3.18}) generate via
the Miura transform (\ref{3.11}) an infinite family of first integrals of 
the mKdV equations.}

For instance, we give the values of the functionals 
(\ref{3.20}) for the unit sphere and the Clifford torus :

1) for the unit sphere :
$$
{\tilde H}_0 = 2, \ \
{\tilde H}_1 = \frac{7}{6}, \ \
{\tilde H}_2 = \frac{2}{5},
$$

2) for the Clifford torus :
$$
{\tilde H}_0 = \pi, \ \
{\tilde H}_1 = \frac{3\sqrt{2}+2}{4}\pi, \ \
{\tilde H}_2 = \frac{\pi}{2}.
$$

\bc
{\bf \S 4. Deformations of tori of revolution via the mKdV flows} 
\ec

We see that for $q = 4U$ and $\lambda=-1$
the operator $L(\lambda)$ (\ref{3.1}) is the same as the operator 
${\tilde L}$ (\ref{2.9}) coming in the general representation formulae for 
surfaces of revolution. The latter one is the reduction of (\ref{2.12}) for
$U(z,{\bar z}) = U(x)$. 

The operator (\ref{2.12}) enters the $L,A,B$-triple representation for the
modified No\-vi\-kov--Veselov (mNV) equation. Since the Novikov--Veselov
(NV) equation is a two-dimensional generalization of  the KdV equation
(\cite{VN}), it is reasonable to expect that there exists a two-dimensional  
generalization of the mKdV equation which is also a modification of
the NV equation in the same sense as the mKdV equation is a modification of 
the KdV equation. This generalization, the mNV equation, was obtained 
by Bogdanov (\cite{Bogd}).  

It was Konopelchenko (\cite{Kon}) who 
observed that if a surface is defined
locally via the generalised Weierstrass representation
from the operator
(\ref{2.12}) and we take the mNV deformation of $U$, 
we obtain a local deformation of a surface. 

Now it is clear that the generalised Weierstrass representation gives
locally all smooth surfaces in ${\bf R^3}$ and 
this result, in fact, goes back to the
19th century and its derivation is almost equivalent to the definition of
the second fundamental form (\cite{T}).
Here we give a strong derivation of this 
fact for surfaces of revolution (Proposition 1). 

We proposed to study these deformations as global deformations of surfaces
and had shown that tori of revolution are deformed into tori 
via the deformation induced by (\ref{3.7}). 
Now this is proved for all mKdV deformations and,
moreover, for the mNV deformation (\cite{T}). Here we describe the global
deformations of tori of revolution via the mKdV hierarchy (Theorem 1).

Explain the definition of the mKdV deformations.

Since the $n$-th mKdV equation is the compatibility condition 
(\ref{3.5}) for (\ref{3.1}) and (\ref{3.2}) the deformation of the
potential $q = 4U$  via (\ref{3.6}) induces the deformation of 
${\tilde \psi}$ via (\ref{3.2}). 
Substituting the deformed function $\psi(x,y,t)$ (\ref{2.5}) into 
(\ref{2.13}) we obtain the {\it $n$-th mKdV deformation of a surface of
revolution}.

{\bf Theorem 1.}
{\sl For every $n \geq 1$ the $n$-th mKdV deformation transforms 
tori of revolution into tori of revolution with preserving 
their conformal types.}   

Proof of Theorem 1.

All tori of revolution are conformally equivalent to rectangular tori
${\bf C}/\{T{\bf Z} + 2\pi\sqrt{-1}{\bf Z}\}$. In terms of (\ref{2.13}),
$T$ is the period of $U(x)$. Since the periods of potentials $U$ are 
preserved 
by the mKdV equations, for proving that the conformal types of tori are 
preserved it suffices to prove that tori are preserved but not closed up
by these deformations. 

By (\ref{2.16}), a preserving of a torus under the $n$-th mKdV deformation is 
equivalent to the following identity
\be
{\cal I}_n = \frac{\partial}{\partial t_n} \int_0^T r(x)s(x) dx \equiv 0.
\label{4.1}
\ee

By (\ref{2.5}), (\ref{3.2}), and (\ref{3.3}), we have
$$
{\cal I}_n = \int_0^T \left(
\frac{\partial r(x)}{\partial t_n}s(x) +
r(x)\frac{\partial s(x)}{\partial t_n} \right) dx   
=
$$
\be
\int_0^T \left(
\sum_0^{n-1} S^{(n)}_{2k+1} \lambda^{2k+1} (r^2 + s^2) +
\sum_0^n T^{(n)}_{2k} \lambda^{2k}(s^2-r^2) \right) dx.
\label{4.2}
\ee

It follows from (\ref{2.7}) that
\be
(rs)_x = -2U(r^2 - s^2), \ \ 
(r^2 + s^2)_x = \lambda (r^2 - s^2), 
\label{4.3}
\ee
$$
(r^2 - s^2)_x = \lambda(r^2+s^2) + 8Urs.
$$
Integrating (\ref{4.2}) by parts, with the use of (\ref{4.3}),
we obtain
\be
{\cal I}_n = \int_0^T = 
\left(
\sum_0^{n-1} \lambda^{2k+1} D^{n-k-1}q_x + 
\sum_0^n \lambda^{2k-1} D^{n-k} q_x \right) (r^2 + s^2) dx
\label{4.4}
\ee
where $q = 4U$. 

Thus, we are left to prove that for every $k \geq 0$
\be
{\cal J}_k = \int_0^T (D^k q_x)(r^2+s^2) dx = 0.
\label{4.5}
\ee
The validity of (\ref{4.1}) for every $n$ follows from (\ref{4.4}) and
(\ref{4.5}) immediately.

From (\ref{4.3}) we derive
\be
{\cal J}_0 = \int_0^T q_x(r^2+s^2) dx  = 
- \int_0^T \lambda q(r^2-s^2) dx =
2\lambda \int_0^T (rs)_x dx = 0.
\label{4.6}
\ee

By (\ref{3.4}), we compute
$$
{\cal J}_k = \int_0^T (D^k q_x)(r^2+s^2) dx = 
$$
\be
\int_0^T(\partial_x^2 + q^2 +q_x\partial_x^{-1}q)(D^{k-1}q_x)(r^2+s^2)dx
= F_1 + F_2 + F_3
\label{4.7}
\ee
where, by virtue of (\ref{4.3}),
$$
F_1 = \int_0^T (D^{k-1}q_x)(r^2+s^2)_{xx} dx = 
$$
$$
\lambda^2 \int_0^T(D^{k-1}q_x)(r^2+s^2)dx + 
2\lambda\int_0^T qrs(D^{k-1}q_x)dx,
$$
$$
F_2 = \int_0^Tq^2(D^{k-1}q_x)(r^2+s^2) dx,
$$
and
$$
F_3 = \int_0^T q_x\partial_x^{-1}(qD^{k-1}q_x)(r^2+s^2) dx =
\int_0^T \partial_x^{-1}(qD^{k-1}q_x)(r^2+s^2) dq = 
$$
$$
- \int_0^T q^2(D^{k-1}q_x)(r^2+s^2) dx - 
\lambda\int_0^T q(r^2-s^2)\partial_x^{-1}(qD^{k-1}q_x) dx =
$$
$$
- \int_0^T q^2(D^{k-1}q_x)(r^2+s^2) dx - 2\lambda\int_0^T qrs(D^{k-1}q_x) dx.
$$
Combining these formulae together we derive
\be
{\cal J}_k = \lambda^2 {\cal J}_{k-1}
\label{4.8}
\ee
and, by (\ref{4.6}), we conclude that the equality
(\ref{4.1}) holds for every $n$.

This proves the theorem.
\footnote[3]
{When we proved this fact first for the mKdV equation (\ref{3.7})
simplifying the problem by reducing the operator (\ref{2.12}) to
(\ref{2.9}), Konopelchenko pointed out to us that the latter operator
enters the $L,A$-pair for the sine-Gordon equation and, thus, 
this equation also induces 
a deformation of surfaces of revolution which is not
generalised yet for all surfaces. However, as we checked, this deformation
closes up tori into cylinders.} 

Since tori of revolution are transformed into tori, the values of 
the functionals ({\ref{3.19})  are correctly defined as the integrals
of forms over closed surfaces and it follows from Proposition 3 that 
these functionals  are preserved by the mKdV deformations.

The first functional ${\tilde H}_0$ is proportional to the Willmore 
functional (\cite{T}), the squared mean curvature, and 
we first derived its conservation (also for the mNV equation
straightforwardly as follows.

Notice that the operator
\be
D^+ = \partial^2_x + q^2 -q \partial^{-1}_x q_x 
\label{4.9}
\ee
is formally coadjoint to $D$, i.e.,
\be
\int f \cdot Dg \ dx = \int D^+f \cdot g \ dx.
\label{4.10}
\ee
and as one can see
$$
\partial_x D^+ = D \partial_x.
$$

Put
$$
W = \int_0^T q^2 dx.
$$
By straightforward computations we derive
$$
\frac{\partial W}{\partial t_n} = 2 \int_0^T q(D^nq_x)dx = 
2 \int_0^T q_x((D^+)^n q) dx = 
$$
$$
2(q(D^+)^n q)|_0^T - 2\int_0^T q\partial_x((D^+)^nq) dx = 
2(q(D^+)^n q)|_0^T - 2\int_0^T q(D^nq_x) dx.
$$
Since the function $q(D^+)^kq$ is periodic, we have
$$
\int_0^T q(D^n q_x) dx = \frac{1}{2} (q(D^+)^n q)|_0^T = 0. 
$$
Thus we conclude that

{\bf Proposition 4.}
{\sl The squared mean curvature $W$ of a torus of revolution is preserved by
the mKdV deformations.}

Describe  the set of stationary points of the mKdV deformation
(\ref{3.7}). Notice that a translation of the argument
$$
\frac{\partial q}{\partial t} = \mbox{const} \cdot 
\frac{\partial q}{\partial x}
$$
also preserves a surface as a geometric object. Hence, we conclude that

{\bf Proposition 5.}
{\sl A surface of revolution corresponding to a potential
$q(x) = 4 U(x)$ via Theorem 1 is preserved by the mKdV deformation
(\ref{3.7}) if and only if
\be
q^2_x = -\frac{q^4}{4} + aq^2 + bq + c
\label{4.11}
\ee
with $a, b,$ and $c$ constants. In this case the deformation reduces to}
$$
q(x,t) = q(x + at).
$$

By (\ref{2.21}), the potential of the unit sphere, $q(x) = 2/\cosh x$
satisfies the equality
\be 
q_x^2 = -\frac{q^4}{4} + q^2,
\label{4.12}
\ee 
and, by (\ref{2.22}), the potential of the round torus $T^2_R$
satisfies the equality
\be
q_x^2 = -\frac{q^4}{4} + (1+\frac{R^2}{2})q^2 + 2Rq + (R^2 - \frac{R^4}{4}).
\label{4.13}
\ee

Since (\ref{4.12}) and (\ref{4.13}) are of the form (\ref{4.11}),
we infer that the unit sphere and the round tori are stationary under the
mKdV deformation (\ref{3.7}). Together with (\ref{3.6}) this implies that

{\bf Proposition 6.}
{\sl The unit sphere and the round tori $T^2_R$ are stationary
points of the mKdV deformations (\ref{3.6}).}

\bc
{\bf \S 5. Geometric transformations of surfaces and 
induced transformations of their potentials}
\ec

It is easy to notice that the homotheties 
$$
{\vec Z} \rightarrow \mbox{const} \cdot {\vec Z}
$$
preserve the potentials of surfaces. 

For surfaces of revolution there are two classes of 
natural transformations else :

1) inversion preserving a structure of a surface of revolution ;

2) a passing to a dual surface in the sense of isothermic surfaces.

The first one takes the form
\be
\Sigma \rightarrow \Sigma_p : \ \
{\vec Z} \rightarrow \frac{{\vec Z}-{\vec P}}{|{\vec Z}-{\vec P}|^2}
\label{5.1}
\ee
where the point $P = (0,0,p)$ lies on the axis of symmetry. 
It follows from (\ref{2.1}) that the image of this transform is also a 
surface of revolution given by (\ref{2.1}) with
\be
\theta(x,p) = \frac{\theta(x)}{\theta(x)^2 + (\varphi(x)-p)^2}, \ \
\varphi(x,p) = \frac{\varphi(x)-p}{\theta(x)^2 +(\varphi(x)-p)^2}.
\label{5.2}
\ee

The induced transform of the potential $q(x)$ of the surface looks more 
difficult. 

First, notice that the formula (\ref{2.4}) can be written in an equivalent
form 
\be
U(x) = \frac{1}{4\theta^2}
(\theta\varphi_x + \theta_x\varphi_{xx}-\theta_{xx}\varphi_x).
\label{5.3}
\ee
This equivalence follows by straightforward computations from
(\ref{2.2}) and (\ref{2.4}).
 
{\bf Proposition 7.}
{\sl The inversion (\ref{5.1}) with $p=0$ transforms the potential $U(x)$
(\ref{5.3}) into the following}
\be
{\tilde U}(x) \ \ ( = {\tilde U}(x,p)|_{p=0}) \ 
= \frac{1}{4\theta^2}(\theta\varphi_x +\theta_{xx}\varphi_x - 
\theta_x\varphi_{xx}) + 
\label{5.4}
\ee
$$
\frac{1}{2\theta^2(\theta^2+\varphi^2)}
((\theta_x^2+\varphi_x^2)(\theta\varphi_x - \theta_x\varphi) -
\theta\varphi(\theta\theta_x + \varphi\varphi_x)).
$$

This transform is an involution but it is not 
represented by a formula in $U(x)$ 
and its derivatives only. For instance, the cylinder
$U(x) = 1/4$ is transformed into the surface with
$$
{\tilde U}(x) = \frac{1}{4} + \frac{1-x^2}{2(1+x^2)}.
$$

We do not give an explicit formula for a generic $p$ but we 
give the following proposition which 
is obtained from (\ref{2.4}) and (\ref{5.2}) by straightforward 
computations.

{\bf Proposition 8.} 
$$
{\tilde U}(x,p) \rightarrow U(x) \ \mbox{as} \ p \rightarrow \infty.
$$

The dual surface $\Sigma^*$ 
is defined for every isothermic surface $\Sigma$, i.e., a surface
with diagonal first and second fundamental forms :
\be
I = e^{2u}( dx^2 + dy^2) , \ \ 
II = e^{2u}( k_1 dx^2 + k_2 dy^2).
\label{5.5}
\ee
The fundamental forms of the dual surface are
\be
I = e^{-2u}( dx^2 + dy^2), \ \
II = -k_1 dx^2 + k_2 dy^2.
\label{5.6}
\ee
The dual surface $\Sigma^* = \{{\vec Z}^*(x,y)\}$
to a surface $\Sigma = \{{\vec Z}(x,y)\}$ is defined up to translations by 
\be
{\vec Z}^*_x = e^{-2u} {\vec Z}_x, \ \
{\vec Z}^*_y = - e^{-2u} {\vec Z}_y.
\label{5.7}
\ee

Every surface of revolution is isothermic and in terms of
$\theta(x)$ the transform to a dual surface looks as
\be
\theta(x) \rightarrow \theta^*(x) = \frac{1}{\theta(x)}.
\label{5.8}
\ee
In this event the function $\varphi^*(x)$ is defined by the formula
\be
\frac{d\varphi^*(x)}{dx} = 
\frac{1}{\theta(x)^2}\frac{d\varphi(x)}{dx}.
\label{5.9}
\ee
 
By (\ref{5.3}), (\ref{5.8}), and (\ref{5.9}), we have

{\bf Proposition 9.}
{\sl The potential $U^*(x)$ of the dual surface
is written as}
\be
U^*(x) = 
\frac{1}{4\theta^2}(\theta\varphi_x - \theta_x\varphi_{xx} + 
\theta_{xx}\varphi_x).
\label{5.10}
\ee 

The formula (\ref{5.10}) is interesting because it appears in
\cite{BHPP} as the square root of the conformal factor of the metric
induced by the central sphere congruence of a surface of revolution
(see the formula (33) in \cite{BHPP} where in our notations 
$k$ equals $2U^*$). 

This shows that there exists an interesting  
relation of the potentials of surfaces to
sphere congruences.

\bc
{\bf \S 6. The conjecture on conformal invariants.}
\ec

Here we would like to introduce the following 

{\bf Conjecture.}
{\sl The integrals (\ref{3.19}) are invariant under the deformations
(\ref{5.1}) of surfaces.}

This means that these integrals are invariant under all deformations
induced by conformal deformations of ${\bf R}^3$ and transforming surfaces 
of revolution into surfaces of revolution.

For ${\tilde H}_0$ this is the Willmore theorem (\cite{W}) on the
invariance of the squared mean curvature (\cite{T}).

We mentioned about this conjecture 
in \cite{T} but for the invariants of general
surfaces, i.e., the first integrals of the mNV flows. The integrals 
(\ref{3.19}) are reductions of them and we may write them in a compact 
form and relate to the well-studied spectral problems. 

If this conjecture is true this means that all the potentials
$q_p$ of the surfaces $\Sigma_p$ generate the isospectral operators   
(\ref{2.9}) (this is deduced from the trace formulae for such operators and 
is usual for the one-dimensional soliton theory).
In virtue of Proposition 8 this implies that we have a circle ${\bf R} \cup
\infty$ of isospectral operators and this circle has to lie on the 
Abelian variety (probably, infinite-dimensional) generated by the mKDV flows.
This strongly relates the M\"obius transform (\ref{5.1}) with 
the mKDV hierarchy.

Since it is evident that this conjecture holds for round tori, we did some 
computations for simple ellipses and obtained the following results :
\be
1) \ \ \ \ \  \frac{x^2}{2}+(y-2)^2 = 1,
\label{6.1}
\ee
$$
4{\tilde H}_0 = 14.733,
$$
$$
\begin{array}{cccccc}
q_p & p=0 & p=1 & p=2 & p=3 & p = \infty \\
& & & & & \\
16{\tilde H}_1 & -31.1181 & -31.1181 & -31.1181 & -31.1181 & -31.1181 \\
32{\tilde H}_2 & 3838.92 & 3838.81 & 3838.81 & 3838.81 & 3838.81
\end{array} ;
$$

\be
2) \ \ \ \ \ \frac{x^2}{3} + (y-2)^2 = 1,
\label{6.2}
\ee
$$
4{\tilde H}_0 = 16.1379, \ \ 32{\tilde H}_2 = 14590.7,
$$
$$
16{\tilde H}_1(q_p) = - 142.454 \ \ \mbox{for} \ \ p = 0,1,2,3,\infty.
$$

These computation are done by Mathematica which does not manage to
proceed with checking the conjecture for ${\tilde H}_2$ in the second case.

Nevertheless, these computations confirm the conjecture.    

This conjecture can be generalised for all surfaces as we mentioned in 
\cite{T}. Now we pass to the problem which concerns isothermic
surfaces only.

Since surfaces of revolution are isothermic, it is interesting to see
how these integrals are transformed by a passing to the dual surface.
Of course, generically the dual surface to a torus of revolution 
is not closed, the group ${\bf Z}$ acts on it by translations, and 
in this case we have to mean by ${\tilde H}_k$ the 
integral over the fundamental domain of this action.

Since the Euler characteristic of a torus vanishes and by 
(\ref{5.5}) and (\ref{5.6}), we have
$$
\int_0^T q^{*2} dx = \int_0^T q^2 dx
$$
for $q^* = 4U^*$ and $q=4U$.
It is more curious that numeric computations by Mathematica shows that
for tori (\ref{6.1}) and (\ref{6.2}) the functionals
${\tilde H}_1$ and ${\tilde H}_2$ are also preserved
(in fact, for the dual surface to the torus (\ref{6.2}) Mathematica gives
${\tilde H}_2 = 14589.9$). The same  is valid for the Clifford torus.

However for the unit sphere we have $\theta^*(x) = \cosh x, \varphi^*(x) = x$,
and $U^*(x) = 0$.

Hence, we pose

{\bf Problem.}
{\sl Describe the class of surfaces of revolution 
for which the well-defined functionals ${\tilde H}_j$ 
are invariant under a
passing to the dual surface.}
 
Probably these are surfaces which do not intersect the axis of revolution,
i.e., with periodic potentials.
In this case these invariants have to be treated as invariants of contours
in the upper-half plane.

We would like to finish with the following piece of speculations.

This conjecture of conformal invariance 
(i.e., invariance with respect to conformal changes of the metric of
the ambient space)
of ${\tilde H}_j$ is a part of
the main conjecture which speaks about conformal invariance of 
the first integrals of higher mNV equations. The latter integrals are
defined for generic immersed surfaces but not only for 
surfaces of revolution. 

These integrals give us new knot invariants as follows. 
Let $\gamma$ be a knot in ${\bf R}^3$. Consider a family 
$M_{\gamma}$ of embedded tori such that every torus $T \in M_{\gamma}$ 
bounds a handlebody $N$ which is represented by an embedding
$f: S^1 \times \{(x,y) \in {\bf R}^2 : x^2 + y^2 \leq 1\} \rightarrow 
{\bf R}^3$ with $f(S^1 \times (0,0)) = \gamma$.
Now put
$$
c_n(\gamma) = \inf_{\Sigma \in M_{\gamma}} |{\tilde H}_n(\Sigma)|.
$$
A very interesting lower estimate for $c_0$ 
(for the case of the Willmore functional) was derived by Willmore
from the unpublished paper of Kearton (\cite{W}).
In particular, the inequality 
$$
c_0(\gamma) \geq 16\pi n
$$
holds where $n$ is the bridge number of $\gamma$.   
It is expected that the higher invariants look like analoguaes 
of hyperbolic volumes. 

{\bf Acknowledgements.} 

The present article was written during author's stay at the Institute
of Theoretical Physics of Freie-Universit\"at in Berlin and was
supported by Volkswagenwerkstiftung and the Russian Foundation for 
Basic Researches (grant 96-01-01889).

The author thanks P. Grinevich, U. Pinkall, J. Richter, and M. Schmidt
for helpful conversations.

\vskip5mm

Institute of Mathematics,

630090 Novosibirsk, Russia

e-mail : taimanov@math.nsc.ru


\begin{thebibliography}{999}

\bibitem{Bogd}
L. V. Bogdanov,
Veselov--Novikov equation as a natural two-dimensional generalization
of the Korteweg--de Vries equation,
Theoret. and Math. Phys. {\bf 70} (1987), 309--314.

\bibitem{BHPP}
F. Burstall, U. Hertrich--Jeromin, F. Pedit, and U. Pinkall,
Curved flats and isothermic surfaces, Preprint, November 1994
(dg-ga 9411010), to appear in Math. Z.

\bibitem{KPP}
G. Kamberov, F. Pedit, and U. Pinkall,
Bonnet pairs and isothermic surfaces, Preprint, October 1996
(dg-ga 9610006).

\bibitem{Kon}
B. G. Konopelchenko,
Induced surfaces and their integrable dynamics,
Studies in Appl. Math. {\bf 96} (1996), 9--52.

\bibitem{KS}
R. Kusner and N. Schmitt,
The spinor representation of surfaces in space,
Preprint, October 1996 (dg-ga 9610005).

\bibitem{New}
A. C. Newell,
Solitons in Mathematics and Physics,
Regional Conference Series in Applied Mathematics {\sl 48},
SIAM, 1985.

\bibitem{Schief}
W. K. Schief, An infinite hierarchy of symmetries associated with
hyperbolic surfaces,
Nonlinearity {\bf 8} (1995), 1--9.

\bibitem{T}
I. A. Taimanov,
Modified Novikov--Veselov equation and differential geometry of surfaces,
Preprint, November 1995 (dg-ga 9511005), 
to appear in Translations of the Amer. Math. Soc.

\bibitem{VN}
A. P. Veselov and S. P. Novikov,
Finite-zone, two-dimensional potential Schr\"odinger operators.
Explicit formulas and evolution equations,
Soviet Math. Dokl. {\bf 30} (1984), 588--591.

\bibitem{W}
T. J. Willmore,
Total Curvature in Riemannian Geometry,
John Wiley and Sons, New York, 1982. 

\end{thebibliography}
\end{document}